\documentclass[12pt]{amsart}%article}%
\usepackage{amsfonts}
\usepackage{amsmath,amssymb,amsthm,amsxtra}
\usepackage{float}
\usepackage[colorlinks, linkcolor=blue,anchorcolor=Periwinkle,
    citecolor=red,urlcolor=Emerald]{hyperref}
\usepackage[usenames,dvipsnames]{xcolor}
\usepackage{enumitem}
\usepackage{geometry,array} 
\usepackage{graphicx}
\usepackage{subfigure}
\usepackage{bookmark}
\usepackage{tikz}\usetikzlibrary{matrix}
\usepackage{url}

\usetikzlibrary{decorations.markings}
\tikzset{->-/.style={decoration={  markings,  mark=at position #1 with
    {\arrow{>}}},postaction={decorate}}}
\tikzset{-<-/.style={decoration={  markings,  mark=at position #1 with
    {\arrow{<}}},postaction={decorate}}}
\usepackage[all]{xy}
\usetikzlibrary{arrows}
\usetikzlibrary{decorations.pathreplacing,decorations.pathmorphing,shadings,fadings,calc}
\usepackage{extarrows}
\def\dexc{cyan!21}%{red!15!blue!40!green}

\theoremstyle{plain}
\newtheorem{theorem}{Theorem}

\theoremstyle{definition}
\newtheorem{definition}[theorem]{Definition}

\theoremstyle{remark}
\newtheorem{remark}[theorem]{Remark}

\numberwithin{equation}{section}

\def\hua{\mathcal}
\def\kong{\mathbb}
\def\<{\langle}
\def\>{\rangle}
\def\ZZ{\mathbb{Z}}

\def\Aut{\operatorname{Aut}}

\def\Hom{\operatorname{Hom}}

\def\Ext{\operatorname{Ext}}

\def\Stab{\operatorname{Stab}}
\def\Stap{\operatorname{Stab}^\circ}

\def\rank{\operatorname{rank}}
\def\numbers{\begin{enumerate}[label=\arabic*{$^\circ$}.]}
\def\ends{\end{enumerate}}

\newcommand{\EG}{\operatorname{EG}}       %exchange graph of heart (oriented)
\newcommand{\EGp}{\operatorname{EG}^\circ}       %principal component
\newcommand{\C}{\hua{C}}

\newcommand{\CEG}[2]{\operatorname{CEG}_{#1}(#2)}             %cluster tilting sets
\newcommand{\D}{\operatorname{\hua{D}}}

\newcommand{\per}{\operatorname{per}}
\newcommand{\Tri}{\bigtriangleup}
\def\arrow{red}
\def\surf{\mathbf{S}}                       %FST's surface
\newcommand{\ST}{\operatorname{ST}}        %spherical twists
\newcommand{\BT}{\operatorname{BT}}        %braid twists
\newcommand{\MCG}{\operatorname{MCG}}

\def\Diff{\operatorname{Diff}}

\def\TT{\kong{T}}
\def\T{T}

\def\M{\mathbf{M}}

\def\surfo{{\mathbf{S}}_\Tri}

\newcommand\Bt[1]{\operatorname{B}_{#1}}

\def\XX{\mathbb{X}}

\def\DMS{\D_{fd}(\surfo)}
\newcommand{\FQuad}{\operatorname{FQuad}}
\newcommand{\numarc}{n}
\newcommand{\numtri}{\aleph}
\newcommand{\uEG}{\underline{\EG}} % unoriented exchange graph of hearts
\newcommand{\SBr}{\operatorname{SBr}}
\newcommand\xx{\mathbf{X}} %compact version of S
\title[Decorated Marked Surfaces]
{Decorated marked surfaces: \\Calabi-Yau categories and related topics}

\author{Yu Qiu}
\address{
\begin{flushleft}
        \hspace{0.3cm}  Ysu Mathematical Sciences Center\\
         \hspace{0.3cm}  Tsinghua University \\
         \hspace{0.3cm}  Beijing, China.\\
\end{flushleft}}
\email{yu.qiu@bath.edu}

\thanks{The survey is in a final form and no version of it will be submitted for publication elsewhere.}

\begin{document}

%%%%%%%%%%%%%%%%%% Abstract Form %%%%%%%%%%%%%%%

\begin{abstract}
This is a survey on the project `Decorated Marked Surfaces'
(\cite{QQ,QQ2,QZ2,BQZ,QZ3,Qs,KQ2}),
where we introduce the decoration $\Delta$ on a marked surfaces $\surf$,
to study Calabi-Yau-2 (cluster) categories, Calabi-Yau-3 (Fukaya) categories,
braid groups for quivers with potential, quadratic differentials and stability conditions.
%The central objects here are the two types of categories associated to $\surfo$,
%namely, the Calabi-Yau-2 (cluster) categories $\hua{C}(\surf)$ and the Calabi-Yau-3 (Fukaya) categories $\DMS$.
%We give a bunch of correspondences between topology and category and
%the punchline of this project is proving the isomorphism $\Stap\DMS\cong\FQuad(\surfo)$
%between complex manifolds and they are simply connected.
%Here $\Stap\DMS$ is the space of stability conditions on $\DMS$ and $\FQuad(\surfo)$
%is the moduli space of $\surfo$-framed quadratic differentials.

% \bigskip

% {\it Key Words:} \quad
%    cluster theory, braid groups, spherical twists, Calabi-Yau categories,
%    stability conditions, quadratic differentials

% \medskip
% {\it $2000$ Mathematics Subject Classification{\rm :}}
% \quad Primary  16Gxx, 16Dxx;  Secondary 16Exx, 16Lxx.
\end{abstract}
%=========================================================
%=========================================================
%\setcounter{tocdepth}{1}
%\tableofcontents\addtocontents{toc}{\setcounter{tocdepth}{1}}
\maketitle
\thispagestyle{empty}
%=========================================================

\def\dgen{2\mathbb{N}_{\leq g}-1}
\def\dgene{2\mathbb{N}_{\leq g}}

%=========================================================
\section{Introduction}
%=========================================================
We introduce the decorated marked surface from Fomin-Shapiro-Thurston's marked surface,
which was originally studied in the cluster theory.
The motivation comes from studying Calabi-Yau categories associated to quivers with potential
and the corresponding spaces of Bridgeland stability conditions.
This project is trying to understand categories via (surface) topology,
which also fits into the framework \cite{DHKK} that relates
dynamical systems (i.e. quadratic differentials) and categories.

%=========================================================
%\subsection*{Context}\label{c}
%=========================================================

\subsection*{Acknowledgements}
This survey is contributed to the proceedings
of the 51st Symposium on Ring Theory and Representation Theory (Sep. 2018).
Qy would like to thank the organizers for the invitation to give lectures in the symposium.

%=========================================================
\section{Categories and topology}\label{sec:MS}
%=========================================================
\subsection{Cluster categories and marked surfaces}
\begin{figure}[ht]\centering
\begin{tikzpicture}[scale=.3]
    \path (-135:4) coordinate (v1)
          (-45:4) coordinate (v2)
          (45:4) coordinate (v3);
\draw[very thick](v1)to(v2)node{$\bullet$}to(v3);
    \path (-135:4) coordinate (v1)
          (45:4) coordinate (v2)
          (135:4) coordinate (v3);
\draw[Emerald, thick](v1)to(v2);
\draw[very thick](v2)node{$\bullet$}to(v3)node{$\bullet$}to
    (v1)node{$\bullet$}(45:1)node[above]{$\gamma$};
\draw[red,thick](135:1.333)node{\tiny{$\circ$}}(-45:1.333)node{\tiny{$\circ$}};
\end{tikzpicture}
\begin{tikzpicture}[scale=1.2, rotate=180]
\draw[blue,<-,>=stealth](3-.6,.7)to(3+.6,.7);
\draw(3,.7)node[below,black]{\footnotesize{in $\surfo$}};
\draw[blue](3-.25,.5-.5)rectangle(3+.25,.5);\draw(3,1.5)node{};
\draw[blue,->,>=stealth](3-.25,.5-.5)to(3+.1,.5-.5);
\draw[blue,->,>=stealth](3+.25,.5)to(3-.1,.5);
\end{tikzpicture}
\begin{tikzpicture}[scale=.3];
    \path (-135:4) coordinate (v1)
          (-45:4) coordinate (v2)
          (45:4) coordinate (v3);
\draw[,very thick](v1)to(v2)node{$\bullet$}to(v3)
(45:1)node[above right]{$\gamma^\sharp$};
    \path (-135:4) coordinate (v1)
          (45:4) coordinate (v2)
          (135:4) coordinate (v3);
\draw[Emerald,,thick](135:4).. controls +(-10:2) and +(45:3) ..(0,0)
                             .. controls +(-135:3) and +(170:2) ..(-45:4);
\draw[,very thick](v2)node{$\bullet$}to(v3)node{$\bullet$}to(v1)node{$\bullet$};
\draw[red,thick](135:1.333)node{\tiny{$\circ$}}(-45:1.333)node{\tiny{$\circ$}};
\end{tikzpicture}

\begin{tikzpicture}[scale=.3]
\draw[thick,>=stealth,->](0,5)to(0,3.6);\draw(0,4.3)node[left]{$^F$};
    \path (-135:4) coordinate (v1)
          (-45:4) coordinate (v2)
          (45:4) coordinate (v3);
\draw[,very thick](v1)to(v2)node{$\bullet$}to(v3);
    \path (-135:4) coordinate (v1)
          (45:4) coordinate (v2)
          (135:4) coordinate (v3);
\draw[,very thick](v2)node{$\bullet$}to(v3)node{$\bullet$}to(v1)node{$\bullet$};
\draw[>=stealth,,thick](-135:4)to(45:4) (45:1)node[above]{$\gamma$};
%\draw[red,thick](135:1.333)node{\tiny{$\circ$}}(-45:1.333)node{\tiny{$\circ$}};
\end{tikzpicture}
\begin{tikzpicture}[scale=1.2, rotate=180]
\draw(3,1.5)node{}(3,.5)node[below]{\footnotesize{in $\surf$}};
\draw[blue,<-,>=stealth](3-.6,.5)to(3+.6,.5);;
\end{tikzpicture}
\begin{tikzpicture}[scale=.3]
\draw[thick,>=stealth,->](0,5)to(0,3.6);\draw(0,4.3)node[right]{$^F$};
    \path (-135:4) coordinate (v1)
          (-45:4) coordinate (v2)
          (45:4) coordinate (v3);
\draw[,very thick](v1)to(v2)node{$\bullet$}to(v3);
    \path (-135:4) coordinate (v1)
          (45:4) coordinate (v2)
          (135:4) coordinate (v3);
\draw[,very thick](v2)node{$\bullet$}to(v3)node{$\bullet$}to(v1)node{$\bullet$};
\draw[>=stealth,,thick](135:4)to(-45:4) (130:1)node[above]{$\gamma^\sharp$};;
%\draw[red,thick](135:1.333)node{\tiny{$\circ$}}(-45:1.333)node{\tiny{$\circ$}};
\end{tikzpicture}
\caption{The decorated/ordinary forward flips}
\label{fig:flip}
\end{figure}
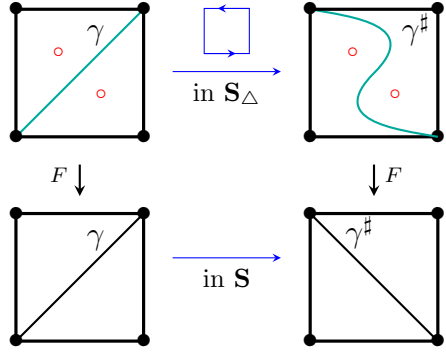

%=========================================================
Following \cite{FST}, a marked surface $\surf$ is a connected oriented smooth surface
with a finite set $\M$ of marked points on its (non-empty) boundary $\partial \surf=\bigcup_{i=1}^b \partial_i$
such that $m_i=\mid \partial_i\cap\M \mid\ge1$.
Such a marked surface is determined, up to diffeomorphism, by the numerical data $g,b$ and the partition
$m=\mid\M\mid=\sum_{i=1}^b m_i$.
We have the following terminology:
\begin{itemize}
\item an \emph{open arc} is (the isotopy class of) a curve on $\surf$ with endpoints in $\M$
but otherwise in $\surf^\circ=\surf \setminus \partial\surf$.
Note that all curves are required to be simple %, i.e. no self-intersection except at the endpoints,
and essential. %, i.e. not homotopic to a constant arc or a boundary arc in $\partial\surf$.
\item two arcs are \emph{compatible} if they do not intersect (except maybe at their endpoints),
\item an \emph{ideal triangulation} $T$ of $\surf$ is a maximal collection of compatible open arcs.
\end{itemize}
An elementary result (cf. \cite[Prop.~2.10]{FST}) is that
any ideal triangulation $T$ of $\surf$, consists of $\numarc=6g-6+3b+m$ open arcs and
divides $\surf$ into $\numtri=({2\numarc+m})/{3}$ triangles.
The \emph{unoriented exchange graph} $\uEG(\surf)$ has vertices
corresponding to ideal triangulations of $\surf$ and edges corresponding to \emph{flips},
as illustrated in the lower row of Figure~\ref{fig:flip}.

Let $\surf$ be a marked surface and $\T$ a triangulation of $\surf$.
Then there is an associated quiver $Q_\T$ with a potential $W_\T$, constructed as follows
(See, e.g. \cite{QZ} for the precise definition):
\begin{itemize}
\item the vertices of $Q_\T$ are (indexed by) the arcs in $\T$;
\item the arrows of $Q_\T$ are (anti-clockwise) angles of triangles in $\T$;
\item these three arrows form a 3-cycle in $Q_\T$ and
$W_\T$ is the sum of all such 3-cycles.
\end{itemize}
The Ginzburg dg algebra $\Gamma_\T=\Gamma(Q_\T, W_\T)$ is the Calabi-Yau-3 algebra
constructed from the combinatorial data.
There are three categories associated to $\Gamma_\T$
\begin{itemize}
\item the \emph{finite dimensional derived category} $\D_{fd}(\Gamma_\T)$, which is Calabi-Yau-3;
\item the \emph{perfect derived category} $\per\Gamma_\T$ (that contains $\D_{fd}(\Gamma_\T)$);
\item the \emph{cluster category} $\C(\Gamma_\T)\colon=\per\Gamma_\T/\D_{fd}(\Gamma_\T)$, which is Calabi-Yau-2.
\end{itemize}

The marked surface $\surf$ provides a good topological model for $\C(\Gamma_\T)$.
In fact, such a category is independent of the choice of $\T$ in the proper sense
and thus will be denoted by $\C(\surf)$.

\begin{theorem}\label{thm:1}\cite{FST,BZ,QZ}
There are following correspondences:
\begin{itemize}
  \item[(O)] The rigid string objects in $\C(\surf)$ are characterized by open arcs on $\surf$.
  \item[(M)] The dimension of $\Ext^1$ between rigid string objects in $\C(\surf)$ is given by
  the intersection numbers between the corresponding open arcs.
  \item[(EG)] The exchange graph $\uEG(\surf)$ is the cluster exchange graph $\CEG{}{\surf}$ for $\C(\surf)$,
  whose the vertices/edges are cluster tilting objects/mutation (cf. \cite{QZ}).
\end{itemize}
\end{theorem}

%=========================================================
\subsection{Calabi-Yau-3 categories and decorated marked surfaces}\label{sec:DMS}
%=========================================================
The first aim of this project is to construct the good topological model for $\D_{fd}(\Gamma_\T)$
and prove the analogue result as Theorem~\ref{thm:1}.
Note that in \cite{BQZ}, we show that $\D_{fd}(\Gamma_\T)$ is independent of $\T$ in a proper sense,
thus we will denote it by $\DMS$.

\begin{figure}[ht]\centering
\begin{tikzpicture}[scale=.2]
  \draw[very thick,NavyBlue](0,0)circle(6)node[above,black]{$_\eta$};
  \draw(-120:5)node{+};
  \draw(-2,0)edge[red, very thick](2,0)  edge[cyan,very thick, dashed](-6,0);
  \draw(2,0)edge[cyan,very thick,dashed](6,0);
  \draw(-2,0)node[white] {$\bullet$} node[red] {$\circ$};
  \draw(2,0)node[white] {$\bullet$} node[red] {$\circ$};
  \draw(0:7.5)edge[very thick,->,>=latex](0:11);\draw(0:9)node[above]{$\Bt{\eta}$};
\end{tikzpicture}\;
%=======================================================
\begin{tikzpicture}[scale=.2]
  \draw[very thick, NavyBlue](0,0)circle(6)node[above,black]{$_\eta$};
  \draw[red, very thick](-2,0)to(2,0);
  \draw[cyan,very thick, dashed](2,0).. controls +(0:2) and +(0:2) ..(0,-2.5)
    .. controls +(180:1.5) and +(0:1.5) ..(-6,0);
  \draw[cyan,very thick,dashed](-2,0).. controls +(180:2) and +(180:2) ..(0,2.5)
    .. controls +(0:1.5) and +(180:1.5) ..(6,0);
  \draw(-2,0)node[white] {$\bullet$} node[red] {$\circ$};
  \draw(2,0)node[white] {$\bullet$} node[red] {$\circ$};
\end{tikzpicture}
\caption{The Braid twist}
\label{fig:Braid twist}
\end{figure}
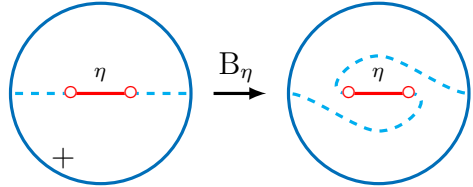

\begin{definition}
The \emph{decorated marked surface} $\surfo$ is a marked surface $\surf$ together with
a fixed set $\Tri$ of $\aleph$ `decorating' points in $\surf^\circ$.
\end{definition}
We also have open arcs in $\surfo$
and a (decorated) triangulation $\TT$ of $\surfo$ is a collection of compatible open arcs
that divides $\surfo$ into $\aleph$ many once-decorated triangles.
The corresponding (oriented) forward flip is shown in the upper row of Figure~\ref{fig:flip},
where one moves the endpoints of an open arc anticlockwise along the quadrilateral to obtain a new open arc
(to form the flipped triangulation).
Denote by $\EG(\surfo)$ the exchange graph of triangulations of $\surfo$,
whose vertices/edges are triangulation of $\surfo$/forward flips.
Note that there is an obvious forgetful map $F\colon\surfo\to\surf$,
which is shown (vertically) in Figure~\ref{fig:flip}.
Then we have the following.
\begin{theorem}\label{thm:2}\cite{QQ2}
There are following correspondences:
\begin{itemize}
  \item[(O)] The reachable rigid string objects in $\per\Gamma_\T$ are characterized by open arcs on $\surfo$.
  \item[(EG)] Any connected component of the exchange graph $\EG(\surfo)$ can be identified
  the principal connected component of the silting exchange graph for $\per\Gamma_\T$,
  where the vertices/edges are silting objects/(forward) mutation (cf. \cite{QQ2}).
\end{itemize}
\end{theorem}
Moreover, we have a new type of arcs:
\begin{itemize}
\item a \emph{closed arc} is (the isotopy class of) a curve on $\surfo$ with different endpoints in $\Tri$.
See red arcs in the left picture of Figure~\ref{fig:Quad A2}.
\end{itemize}
This type of arcs plays a key role in braid groups and mapping class groups of surfaces.
Namely, we have the following notions.
The mapping class group $\MCG(\surfo)$ is the group of isotopy classes of diffeomorphisms of $\surfo$,
where all diffeomorphisms and isotopies fix $\M$ and $\Tri$ setwise.
On the other hand, the mapping class group $\MCG(\surf)$  fixes just $\M$ setwise.
Thus there is a forgetful group homomorphism
\begin{equation}\label{eq:FM}
F_M\colon\MCG(\surfo)\to\MCG(\surf),
\end{equation}
whose kernel is the so-called surface braid group $\SBr(\surfo)$.
We are particularly interested in the braid twist group $\BT(\surfo)$, which is subgroup of $\SBr(\surfo)$,
generated by braid twist along closed arcs, cf. Figure~\ref{fig:Braid twist}.

Then we have the following.
\begin{theorem}\label{thm:3}\cite{QQ,QZ2}
There are following correspondences:
\begin{itemize}
  \item[(O)] The reachable spherical objects in $\DMS$ are characterized by closed arcs on $\DMS$.
  \item[(M)] The dimension of $\Hom^\bullet$ between spherical objects in $\DMS$ is given by
  the intersection numbers between the corresponding closed arcs.
  \item[(G)] The braid twist group $\BT(\surfo)$ can be identified with
  the (so-called spherical twist) subgroup $\ST(\surfo)$ of the auto-equivalence group $\Aut\DMS$,
  where the braid twist along a closed arc becomes the spherical twist of the corresponding spherical object.
  \item[(EG)] Any connected component $\EGp(\surfo)$ of the exchange graph $\EG(\surfo)$ can be identified
  the principal connected component $\EGp\DMS$ of the exchange graph of $\DMS$,
  where the vertices/edges are hearts/simple (forward) tilting (cf. \cite{QQ}).
\end{itemize}
\end{theorem}

%=========================================================
%\subsection{Braid groups and spherical twist groups}\label{sec:Br}
%=========================================================
\begin{remark}\label{rk:2}
The isomorphism $\BT(\surfo)\cong\ST(\surfo)$ mentioned above plays a key role in understanding
the generalization of braid group to the case of quivers with (super)potential, see \cite{Qs} for details.
Their (finite) presentations in \cite{QZ3} also plays a crucial technical role the dissuasion below.
\end{remark}

%=========================================================
\section{Quadratic differentials as stability conditions}\label{sec:Stab=Quad}
%=========================================================
One of the motivations to introduce decorated marked surface is coming from
the interplay between quadratic differentials and categories.
\subsection{Framed uadratic differentials}\label{sec:Quad}
%=========================================================
Let $\xx$ be a compact Riemann surface and $\omega_\xx$
be its holomorphic cotangent bundle.
A \emph{meromorphic quadratic differential} $\phi$ on $\xx$ is a meromorphic section
of the line bundle $\omega_{\xx}^{\otimes2}$.
We consider \emph{GMN differentials} $\phi$ on $\xx$, which are meromorphic quadratic differential
such that all zeroes of $\phi$ are simple and all poles of $\phi$ have order at least three.
The real (oriented) blow-up of $(\xx,\phi)$ is a differentiable surface $\xx^\phi$,
which is obtained from the underlying differentiable surface by replacing each pole $P$ of $\phi$ by a boundary $\partial_P$,
where the points on the boundary correspond to the real tangent directions at $P$.
Furthermore, we will mark the points on $\partial_P$ that correspond to the distinguished tangent directions,
so there are $\operatorname{ord}_\phi(P)+2$ marked points on $\partial_P$.

\begin{figure}[ht]\centering
\begin{tikzpicture}[scale=.5]
    \path (18+72:5) coordinate (v2)
          (18+72*3:5) coordinate (v1)
          (18+72*2:2.7) coordinate (v3)
          (0,-1) coordinate (v4);
  \foreach \j in {.1, .18, .26, .34, .42, .5,.58, .66, .74, .82, .9}
    {
      \path (v3)--(v4) coordinate[pos=\j] (m0);
      \draw[Emerald!60, thin] plot [smooth,tension=.3] coordinates {(v1)(m0)(v2)};
    }
\draw[thick](v4)to(v1)to(v3)to(v2)to(v4);

    \path (18+72:5) coordinate (v2)
          (18+72*4:5) coordinate (v1)
          (18+72*0:2.7) coordinate (v3);
  \foreach \j in {.1, .18, .26, .34, .42, .5,.58, .66, .74, .82, .9}
    {
      \path (v3)--(v4) coordinate[pos=\j] (m0);
      \draw[Emerald!60, thin] plot [smooth,tension=.3] coordinates {(v1)(m0)(v2)};
    }
\draw[thick](v4)to(v1)to(v3)to(v2)to(v4);
%==
    \path (18+72:5) coordinate (v2)
          (18+72*2:5) coordinate (v1)
          (18+72*2:2.7) coordinate (v3);
    \path (v1)--(v2) coordinate[pos=.4] (v4);
  \foreach \j in {.2,.32,.45,.55,.68,.8, .9}
    {
      \path (v3)--(v4) coordinate[pos=\j] (m0);
      \draw[Emerald!60, thin] plot [smooth,tension=.4] coordinates {(v1)(m0)(v2)};
    }
%==
    \path (18+72*3:5) coordinate (v2);
    \path (v1)--(v2) coordinate[pos=.4] (v4);
  \foreach \j in {.2,.32,.45,.55,.68,.8, .9}
    {
      \path (v3)--(v4) coordinate[pos=\j] (m0);
      \draw[Emerald!60, thin] plot [smooth,tension=.4] coordinates {(v1)(m0)(v2)};
    }
%==
    \path (18+72:5) coordinate (v2)
          (18+72*0:5) coordinate (v1)
          (18+72*0:2.7) coordinate (v3);
    \path (v1)--(v2) coordinate[pos=.4] (v4);
  \foreach \j in {.2,.32,.45,.55,.68,.8, .9}
    {
      \path (v3)--(v4) coordinate[pos=\j] (m0);
      \draw[Emerald!60, thin] plot [smooth,tension=.4] coordinates {(v1)(m0)(v2)};
    }
%==
    \path (18+72*4:5) coordinate (v2);
    \path (v1)--(v2) coordinate[pos=.4] (v4);
  \foreach \j in {.2,.32,.45,.55,.68,.8, .9}
    {
      \path (v3)--(v4) coordinate[pos=\j] (m0);
      \draw[Emerald!60, thin] plot [smooth,tension=.4] coordinates {(v1)(m0)(v2)};
    }
%==
    \path (18+72*3:5) coordinate (v2)
          (18+72*4:5) coordinate (v1)
          (0,-1) coordinate (v3);
    \path (v1)--(v2) coordinate[pos=.5] (v4);
  \foreach \j in {.2,.3,.4,.5,.6,.7,.8,.9}
    {
      \path (v3)--(v4) coordinate[pos=\j] (m0);
      \draw[Emerald!60, thin] plot [smooth,tension=.4] coordinates {(v1)(m0)(v2)};
    }
%==
\draw[thick](18+72*2:5)to(18+72*2:2.7)(18:5)to(18:2.7);
\foreach \j in {1,2,3,4,5}
    {\draw[Emerald!50,very thick]
     (18+72*\j:5)to(90+72*\j:5);}
\foreach \j in {1,2,3,4,5}
    {\draw[NavyBlue]
        (18+72*\j:5)node{$\bullet$};}
\draw[red,very thick]
  (18+72*2:2.7)node{$\bullet$}node[white]{\tiny{$\bullet$}}node{\tiny{$\circ$}}to
  (0,-1)node{$\bullet$}node[white]{\tiny{$\bullet$}}node{\tiny{$\circ$}}to
  (18+72*0:2.7)node{$\bullet$}node[white]{\tiny{$\bullet$}}node{\tiny{$\circ$}};
\end{tikzpicture}
\quad
\begin{tikzpicture}[scale=.6,rotate=18]
\draw[dashed, NavyBlue](0,0)circle(3.6)(-90:4.5);
\foreach \j in {0,1,2,3,4}
{
  \draw[very thick,NavyBlue](72*\j+0:3.6)to(0,0);
  \draw[Emerald] plot [smooth,tension=1] coordinates
    {(0,0)(15+72*\j:3)(72-15+72*\j:3)(0,0)};
  \draw[Emerald] plot [smooth,tension=1] coordinates
    {(0,0)(20+72*\j:2.5)(72-20+72*\j:2.5)(0,0)};
  \draw[Emerald] plot [smooth,tension=1] coordinates
    {(0,0)(25+72*\j:2)(72-25+72*\j:2)(0,0)};
  \draw[Emerald] plot [smooth,tension=1] coordinates
    {(0,0)(30+72*\j:1.5)(72-30+72*\j:1.5)(0,0)};
}
\draw[NavyBlue](0,0)node{$\bullet$};
\end{tikzpicture}
\caption{The foliation of a GMN differential on $\kong{P}^1$ for the $A_2$ case}
\label{fig:Quad A2}
\end{figure}
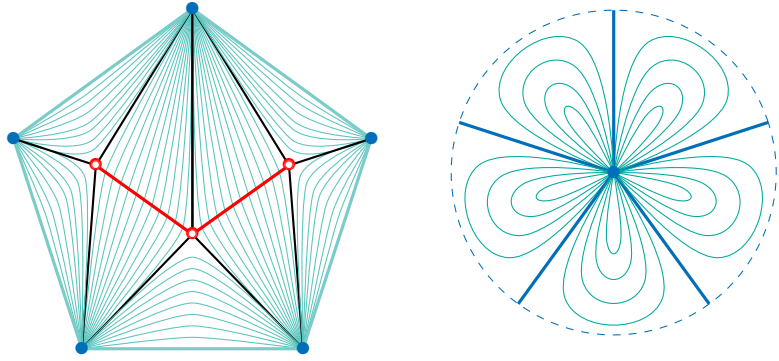

\begin{definition}\label{def:SoFQuad}
The decorated real blow-up $\xx^\phi_\Tri$ of $(\xx,\phi)$ is
the decorated marked surface obtained from $\xx^\phi$ by adding the set of zeroes of $\phi$ as decorations.
Given any decorated marked surface $\surfo$,
an \emph{$\surfo$-framed quadratic differential} $(\xx,\phi,\psi)$
is a Riemann surface $\xx$ with a GMN differential $\phi$,
equipped with a diffeomorphism $\psi\colon\surfo \to\xx^\phi_\Tri$, preserving the marked points and decorations
(see right picture of Figure~\ref{fig:Quad A2})..

Two $\surfo$-framed quadratic differentials $(\xx_1,\phi_1,\psi_1)$ and $(\xx_2,\phi_2,\psi_2)$
are equivalent, if there exists a biholomorphism $f\colon\xx_1\to\xx_2$
such that $f^*(\phi_2)=\phi_1$ and furthermore $\psi_2^{-1}\circ f_*\circ\psi_1\in\Diff_0(\surfo)$,
where $f_*\colon\xx_1^{\phi_1}\to\xx_2^{\phi_2}$ is the induced diffeomorphism.
Here $\Diff_0(\surfo)$ is the identity component of the group
$\Diff(\surfo)$ of diffeomorphisms preserving marked points and decorations (each setwise).
\end{definition}

Denote by $\FQuad(\surfo)$ the moduli space of $\surfo$-framed quadratic differentials
and take any connected component $\FQuad^\circ(\surfo)$ (which are isomorphic to each other).
\begin{remark}\label{rk:1}
The foliation of a quadratic differential gives (the so-called WKB) a triangulation of its decorated real below-up
(see left picture of Figure~\ref{fig:Quad A2}).
Therefore the exchange graph $\EGp(\surfo)$ of triangulations is the skeleton of $\FQuad^\circ(\surfo)$,
as $\ZZ^n$ is the skeleton of $\mathbb{C}^n$ (cf. Figure~\ref{fig:cover1}).
\end{remark}

%=========================================================
\subsection{Stability conditions}\label{sec:Stab}
%=========================================================
\begin{definition}[Bridgeland]
A \emph{stability condition} $\sigma = (Z,\hua{P})$ on a triangulated category $\D$ consists of
a central charge $Z\in\Hom_{\ZZ}(K(\D),\mathbb{C})$ and a collection of
full additive subcategories $\{ \hua{P}(\varphi) \subset \D\mid \varphi\in\kong{R}\}$, such that
\begin{itemize}
\item if $0 \neq E \in \hua{P}(\varphi)$
then $Z(E) = m(E) \exp(\varphi  \pi \mathbf{i} )$ for some $m(E) \in \kong{R}_{>0}$;
\item $\hua{P}(\varphi+1)=\hua{P}(\varphi)[1]$, for all $\varphi \in \kong{R}$;
\item if $\varphi_1>\varphi_2$, then $\Hom(\hua{P}(\varphi_1),\hua{P}(\varphi_2))=0$;
\item any object $E \in \D$ admits a Harder-Narashimhan filtration of triangles (cf. \cite[(8.1)]{QQ}),
whose factors are $A_j \in \hua{P}(\varphi_j)$ for $\leq j\le m$ satisfying
$\varphi_1 > \varphi_2 > ... > \varphi_m$.
\end{itemize}
\end{definition}
Bridgeland shows that all stability conditions on $\D$ form a complex manifold $\Stab\D$ of dimension $n$
with local coordinate $Z$ for $n=\rank K(\D)(<\infty)$.
For $\DMS$, denote by $\Stap\DMS$ its principal connected component.
As in Remark~\ref{rk:1}, exchange graph $\EGp\DMS$ of hearts is the skeleton of $\Stap\DMS$.
Such a philosophy/approach has been explored carefully (cf. the survey \cite{Q5}).
Building on Bridgeland-Smith's seminar work \cite{BS}, we prove the following.
\begin{theorem}\cite{KQ2}
There is an isomorphism $\Stap\DMS\cong\FQuad(\surfo)$ between complex manifolds and they are simply connected.
\end{theorem}

Note that the correspondence (O) in Theorem~\ref{thm:3} plays a key role in the isomorphism above:
the saddle trajectories $\eta$ of a quadratic differential, which are closed arcs,
correspond to semistable(/spherical) objects for stability conditions.

\def\graphc{JungleGreen}
\def\vertexc{blue}
\def\halfc{blue}%{Periwinkle}
\def\thirdc{blue}%{Periwinkle}
\def\dexc{black!10!blue!50!green}
\begin{figure}[ht]\centering
\newcommand{\vsource}{\otimes}
\newcommand{\vsink}{\odot}
\newcommand{\vertx}{\bullet}
\begin{tikzpicture}[scale=1.7,arrow/.style={->,>=stealth}]
\draw[white, fill=\graphc!19] (2.5,2.5) rectangle (3.5,3.5);
\foreach \j in {1.5,2.5,3.5,4.5} {\draw[\graphc!30,thick](\j,1)to(\j,5)(1,\j)to(5,\j);}
\foreach \i in {2,3,4}{
\foreach \j in {1,2,3,4,5}
   {\draw[white] (\i,\j) node (t\j) { };}
\foreach \a/\b in {1/2,2/3,3/4,4/5}{
  \draw[blue, thick] (t\a) edge[arrow] (t\b);
}}
\foreach \i in {2,3,4}{
\foreach \j in {1,2,3,4,5}
   {\draw[white] (\j,\i) node (t\j) { };}
\foreach \a/\b in {1/2,2/3,3/4,4/5}{
  \draw[blue, thick] (t\a) edge[arrow] (t\b);
}}
\foreach \i in {2,3,4}{
\foreach \j in {2,3,4}
   {\draw[blue, thick] (\j,\i) node (t\j) {\scriptsize{$\bullet$}};}}
\end{tikzpicture}
\begin{tikzpicture}[scale=.65]
\path (0,0) coordinate (O);
\path (0,6) coordinate (S1);
\draw[fill=\graphc!19,dotted]
    (S1) arc(360/3-360/3+180:360-360/3:10.3923cm) -- (O) -- cycle;
\draw[,dotted] (210:6cm) -- (O);
\draw[fill=black] (30:1.6077cm) circle (.05cm);
\draw
    (O)[Periwinkle,dotted,thick]   circle (6cm);
\draw
    (S1)[\graphc,thick]
    \foreach \j in {1,...,3}{arc(360/3-\j*360/3+180:360-\j*360/3:10.3923cm)} -- cycle;
\draw
    (S1)[\graphc,semithick]
    \foreach \j in {1,...,6}{arc(360/6-\j*360/6+180:360-\j*360/6:3.4641cm)} -- cycle;
\draw
    (S1)[\graphc]
    \foreach \j in {1,...,12}{arc(360/12-\j*360/12+180:360-\j*360/12:1.6077cm)}
        -- cycle;
\foreach \j in {1,...,3}
{\path (-90+120*\j:.7cm) node[\vertexc] (v\j) {\scriptsize{$\bullet$}};
 \path (-210+120*\j:.7cm) node[\vertexc] (w\j) {\scriptsize{$\bullet$}};
 \path (-90+120*\j:2.2cm) node[\vertexc] (a\j) {\scriptsize{$\bullet$}};
 \path (-90+15+120*\j:3cm) node[\vertexc] (b\j) {\scriptsize{$\bullet$}};
 \path (-90-15+120*\j:3cm) node[\vertexc] (c\j) {\scriptsize{$\bullet$}};}
\foreach \j in {1,...,3}
{\path[->,>=latex] (v\j) edge[\thirdc,bend left,thick] (w\j);
 \path[->,>=latex] (a\j) edge[\thirdc,bend left,thick] (b\j);
 \path[->,>=latex] (b\j) edge[\thirdc,bend left,thick] (c\j);
 \path[->,>=latex] (c\j) edge[\thirdc,bend left,thick] (a\j);}
\foreach \j in {1,...,3}
{\path (60*\j*2-1-60:3.9cm) node[\vertexc] (x1\j) {\scriptsize{$\bullet$}};
% \path (60*\j*2+10:4.4cm) node[\vertexc] (y1\j) {\scriptsize{$\bullet$}};
% \path (60*\j*2-10:4.4cm) node[\vertexc] (z1\j) {\scriptsize{$\bullet$}};
 \path (60*\j*2+1-120:3.9cm) node[\vertexc] (x2\j) {\scriptsize{$\bullet$}};
% \path (60+60*\j*2+10:4.4cm) node[\vertexc] (y2\j) {\scriptsize{$\bullet$}};
% \path (60+60*\j*2-10:4.4cm) node[\vertexc] (z2\j) {\scriptsize{$\bullet$}};
}
\foreach \j in {1,...,3}
{\path[->,>=latex] (v\j) edge[\halfc,bend left,thick] (a\j);
 \path[->,>=latex] (a\j) edge[\halfc,bend left,thick] (v\j);
 \path[->,>=latex] (b\j) edge[\halfc,bend left,thick] (x1\j);
 \path[->,>=latex] (c\j) edge[\halfc,bend left,thick] (x2\j);
 \path[->,>=latex] (x1\j) edge[\halfc,bend left,thick] (b\j);
 \path[->,>=latex] (x2\j) edge[\halfc,bend left,thick] (c\j);}
\end{tikzpicture}
\caption{Skeleton illustrations (for exchange graphs in moduli/stability spaces)}
\label{fig:cover1}
\end{figure}
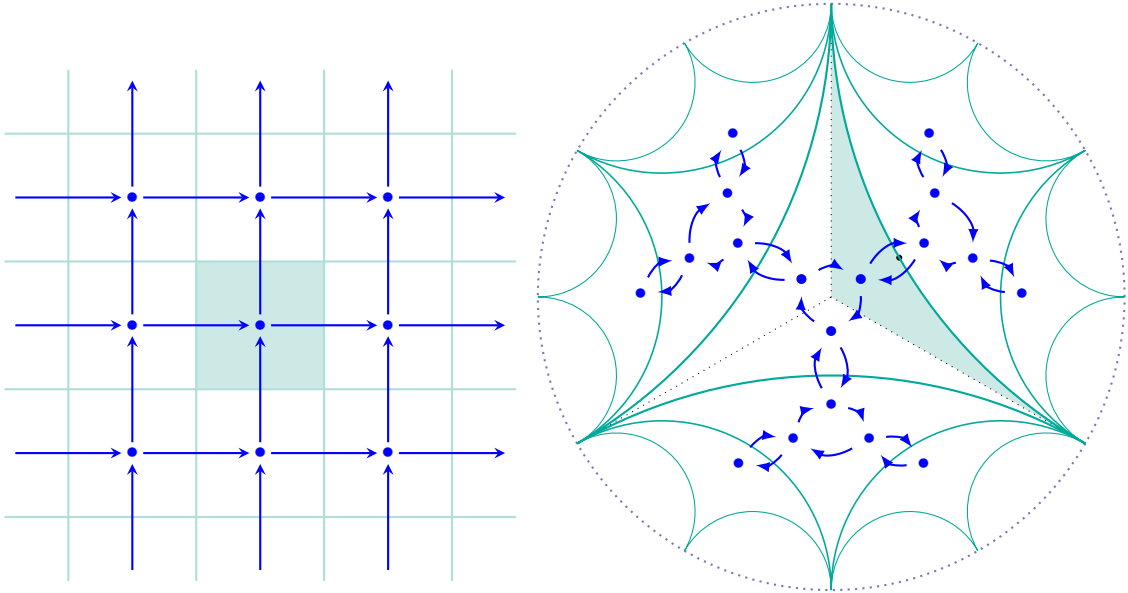

%=========================================================
\section{Further direction}\label{sec:IQ}
%=========================================================
We (Akishi Ikeda and me) have started a new project `$q$-stability conditions on Calabi-Yau $\XX$ categories' (\cite{IQ1,IQ2}),
where we did the following:
\begin{itemize}
  \item In \cite{IQ1}, we set up the framework of $q$-deformation of stability conditions on (Calabi-Yau-)$\XX$ categories
  $\D_\XX$, whose Grotendieck group is $R^{\oplus n}$ for $R=\ZZ[q^{\pm1}]$.
  \item In \cite{IQ2}, we study the surface case and introduce the $q$-deformation of quadratic differentials to
  realize the $q$-deformation of stability conditions in this case.
\end{itemize}
In the forthcoming paper (together with Yu Zhou), we will study the decorated marked surfaces for Calabi-Yau-$\XX$ categories.

%=========================================================
%\clearpage

%=========================================================
\end{document}